\documentclass[lang = american]{ems-icm} %% change to `american' if you use American English
\usepackage[all]{xy}

\theoremstyle{plain}
\newtheorem{theorem}{Theorem}[section]
\newtheorem{prop}[theorem]{Proposition}
\newtheorem{lemma}[theorem]{Lemma}
\newtheorem{corollary}[theorem]{Corollary}

\theoremstyle{definition}
\newtheorem{com}[theorem]{Comment}
\newtheorem{apl}[theorem]{Application}
\newtheorem{exercise}[theorem]{Exercise}
\newtheorem{redu}[theorem]{Reduction}
\newtheorem{refinement}[theorem]{Refinement}
\newtheorem{summary}[theorem]{Summary}
\newtheorem{importnota}[theorem]{Important Notation}
\newtheorem{prblm}[theorem]{Problem}
\newtheorem{notation}[theorem]{Notation}
\newtheorem{explanation}[theorem]{Explanation}
\newtheorem{defin}[theorem]{Definition}
\newtheorem{caution}[theorem]{Caution}
\newtheorem{remark}[theorem]{Remark}
\newtheorem{reminder}[theorem]{Reminder}
\newtheorem{illustration}[theorem]{Illustration}
\newtheorem{observation}[theorem]{Observation}
\newtheorem{construction}[theorem]{Construction}
\newtheorem{discussion}[theorem]{Discussion}
\newtheorem{example}[theorem]{Example}
\newtheorem{motexample}[theorem]{Motivating Example}
\newtheorem{conclusion}[theorem]{Conclusion}
\newtheorem{sketch}[theorem]{Sketch}
\newtheorem{triviality}[theorem]{Triviality}
\newtheorem{proto}[theorem]{Prototype Quasifibration}
\newtheorem{cauex}[theorem]{Cautionary Example}
\newtheorem{hypo}[theorem]{Hypothesis}
\newtheorem{subth}{ }[theorem]
\newtheorem{case}{Case}[theorem]
\newtheorem{ssubth}{ }[subth]
\newtheorem{facts}[theorem]{Facts}
\newtheorem{history}[theorem]{Historical Survey}
\newtheorem{blackbox}[theorem]{Black Box}
\newtheorem{proofs}[theorem]{Discussion of the Proofs, Old and New}

\newcommand\nc{\newcommand}

\nc\tri[1]{\begin{triviality}
\label{#1}}
\nc\bbx[1]{\begin{blackbox}
\label{#1}}
\nc\fac[1]{\begin{facts}
\label{#1}
}
\nc\app[1]{\begin{apl}
\label{#1}
}
\nc\skt[1]{\begin{sketch}
\label{#1}
}
\nc\hst[1]{\begin{history}
\label{#1}
}
\nc\pfs[1]{\begin{proofs}
\label{#1}
}
\nc\cas[1]{\begin{case}
\label{#1}
}
\nc\rfn[1]{\begin{refinement}
\label{#1}}
\nc\prt[1]{\begin{proto}
\label{#1}}
\nc\lem[1]{\begin{lemma}
\label{#1}}
\nc\pro[1]{\begin{prop}
\label{#1}}
\nc\thm[1]{\begin{theorem}
\label{#1}}
\nc\dis[1]{\begin{discussion}
\label{#1}
}
\nc\cor[1]{\begin{corollary}
\label{#1}}
\nc\dfn[1]{\begin{defin}
\label{#1}}
\nc\sthm[1]{\begin{subth}
\label{#1}}
\nc\exm[1]{\begin{example}
\label{#1}
}
\nc\mxm[1]{\begin{motexample}
\label{#1}
}
\nc\obs[1]{\begin{observation}
\label{#1}
}
\nc\plm[1]{\begin{prblm}
\label{#1}
}
\nc\rmk[1]{\begin{remark}
\label{#1}
}
\nc\rmd[1]{\begin{reminder}
\label{#1}
}
\nc\ntn[1]{\begin{notation}
\label{#1}
}
\nc\exe[1]{\begin{exercise}
\label{#1}
}
\nc\xpl[1]{\begin{explanation}
\label{#1}
}
\nc\smr[1]{\begin{summary}
\label{#1}
}
\nc\cau[1]{\begin{caution}
\label{#1}
}
\nc\hyp[1]{\begin{hypo}
\label{#1}}
\nc\imn[1]{\begin{importnota}
\label{#1}
}
\nc\rdn[1]{\begin{redu}
\label{#1}
}
\nc\cax[1]{\begin{cauex}
\label{#1}
}
\nc\cmt[1]{\begin{com}
\label{#1}
}
\nc\con[1]{\begin{construction}
\label{#1}
}
\nc\ill[1]{\begin{illustration}
\label{#1}
}
\nc\ssthm[1]{\begin{ssubth}
\label{#1}
}
\nc\cnc[1]{\begin{conclusion}
\label{#1}
}

\nc\elem{\end{lemma}}
\nc\ebbx{\end{blackbox}}
\nc\erdn{\end{redu}}
\nc\erfn{\end{refinement}}
\nc\eprt{\end{proto}}
\nc\ethm{\end{theorem}}
\nc\ecor{\end{corollary}}
\nc\edfn{\end{defin}}
\nc\esthm{\end{subth}}
\nc\epro{\end{prop}}
\nc\etri{\end{triviality}}
\nc\eexm{
\end{example}}
\nc\emxm{
\end{motexample}}
\nc\eobs{
\end{observation}}
\nc\ecmt{
\end{com}}
\nc\efac{
\end{facts}}
\nc\eapp{
\end{apl}}
\nc\ermk{
\end{remark}}
\nc\ermd{
\end{reminder}}
\nc\eill{
\end{illustration}}
\nc\eplm{
\end{prblm}}
\nc\ecas{
\end{case}}
\nc\eskt{
\end{sketch}}
\nc\ecau{
\end{caution}}
\nc\ecax{
\end{cauex}}
\nc\eimn{
\end{importnota}}
\nc\entn{
\end{notation}}
\nc\eexe{
\end{exercise}}
\nc\expl{
\end{explanation}}
\nc\edis{
\end{discussion}}
\nc\econ{
\end{construction}}
\nc\esmr{
\end{summary}}
\nc\ehst{
\end{history}}
\nc\epfs{
\end{proofs}}
\nc\ehyp{
\end{hypo}}
\nc\ecnc{
\end{conclusion}}
\nc\essthm{
\end{ssubth}}

\nc\script{\mathscr}

\nc\ct{{\script T}}
\nc\cf{{\script F}}
\nc\cg{{\script G}}
\nc\ch{{\script H}}
\nc\ck{{\script K}}
\nc\cl{{\script L}}
\nc\cv{{\script V}}
\nc\cw{{\script W}}
\nc\ce{{\script E}}
\nc\cs{{\script S}}
\nc\car{{\script R}}
\nc\cd{{\script D}}
\nc\cc{{\script C}}
\nc\ca{{\script A}}
\nc\ci{{\script I}}
\nc\cj{{\script J}}
\nc\co{{\script O}}
\nc\cu{{\script U}}
\nc\cx{{\script X}}
\nc\Cp{{\script P}}
\nc\cq{{\script Q}}
\nc\cy{{\script Y}}
\nc\cz{{\script Z}}
\nc\Hom{\mathrm{Hom}}
\nc\la{\longrightarrow}
\nc\zz{\mathbb{Z}}
\nc\nn{\mathbb{N}}
\nc\sS{\mathbb{S}}
\nc\cC{\mathbb{C}}
\nc\rr{\mathbb{R}}
\nc\genu[3]{{\langle#1\rangle}^{[#3]}_{#2}}
\nc\genuf[2]{{\langle#1\rangle}^{}_{#2}}
\nc\ogenu[3]{\ov{\langle#1\rangle}^{[#3]}_{#2}}
\nc\ogenuf[2]{\ov{\langle#1\rangle}^{}_{#2}}
\nc\ov{\overline}
\nc\scl{\text{SmCl}}
\nc\lcl{\text{LgCl}}
\nc\dperf[1]{\D^{\mathbf{perf}}(#1)}
\nc\dcoh{\mathbf{D}^b_{\mathbf{coh}}}
\nc\dcohm{\mathbf{D}^-_{\mathbf{coh}}}
\newcommand{\Dqc}{{\mathbf D_{\mathbf{qc}}}}
\newcommand{\Dqcmi}{{\mathbf D_{\mathbf{qc}}^-}}
\newcommand{\Dqcpl}{{\mathbf D_{\mathbf{qc}}^+}}
\newcommand{\Dqcb}{{\mathbf D_{\mathbf{qc}}^b}}

\nc\dmcoh{\mathbf{D}^-_{\mathrm{coh}}}
\nc\dscoh{\mathbf{D}^{}_{\mathrm{coh}}}
\nc\mmod[1]{#1\text{--}\mathrm{mod}}
\nc\Mod[1]{#1\text{--}\mathrm{Mod}}
\nc\MMod[1]{\text{Mod--}#1}
\nc\D{\mathbf{D}}
\nc\op{^{\mathrm{op}}}
\nc\an{^{\mathrm{an}}}
\nc\spec[1]{\mathrm{Spec}\left(#1\right)}
\nc\prf{\begin{proof}}
\nc\eprf{\end{proof}}
\nc\ogenul[3]{\ov{\langle#1\rangle}^{(-\infty,#3]}_{#2}}
\nc\be{\begin{enumerate}}
\nc\ee{\end{enumerate}}
\nc\wt{\widetilde}
\nc\wi{\wt{\text{\i}}}
\nc\ab{\script{Ab}}
\nc\fl{\mathfrak{L}}
\nc\fs{\mathfrak{S}}
\nc\e{\varepsilon}
\nc\wh{\widehat}

\numberwithin{equation}{section}

\begin{document}

\volumetitle{ICM 2022} % Don't alter this line.

%------
% Insert the title of your paper and (if necessary)
% a short title for the running head.
%------
\title{Finite approximations as a tool for studying triangulated categories}
%\titlemark{}

%------
% Insert full names of the authors.
% Add further authors as follows:
%  \emsauthor{2}{NAME INCL. FULL FIRST NAME}{NAME WITH FIRST NAME INITIALS}
%  \emsauthor{3}{NAME INCL. FULL FIRST NAME}{NAME WITH FIRST NAME INITIALS}
% etc.
%------
\emsauthor{1}{Amnon Neeman}{A.~Neeman}

%------
% Add one \emsaffil and one \email for each author.
%------
\emsaffil{1}{Mathematical Sciences Institute,
        Building 145,
        The Australian National University,
        Canberra, ACT 2601,
        AUSTRALIA \email{Amnon.Neeman@anu.edu.au}}

%\dedication{Dedicated to ...}

%------
% Insert your abstract.
%------
\begin{abstract}
Small, finite entities are easier and simpler to
manipulate than gigantic, infinite
ones. Consequently
huge chunks of mathematics are devoted to methods 
reducing the study of big,
cumbersome
objects to an analysis of their finite building blocks. 
The manifestation of this general pattern, in the study
of derived and triangulated categories, dates  back almost to 
the beginnings of the subject---more precisely to
articles by
Illusie in SGA6, way back in the early 1970s.

What's new, at least new in the world of derived and triangulated
categories, is that one gets extra mileage from analysing more carefully and
quantifying more precisely just
how efficiently one can estimate infinite objects by finite ones.
This leads one to the study of metrics on triangulated categories, and
of how accurately an object can be approximated 
by finite objects of bounded size. 
\end{abstract}

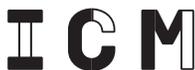
\maketitle

%------
% INSERT THE BODY OF THE PAPER HERE (except
% acknowledgments, funding info and bibliography)
%------

\section{Introduction}
\label{S1}

In every branch of mathematics we try to solve complicated problems
by reducing to simpler ones, and from antiquity people have used
finite approximations to study infinite objects. Naturally, whenever
a new field comes into being, one of the 
first developments is to try to understand what should
be the right notion of finiteness in the discipline.
Derived and triangulated
categories were introduced by Verdier in his PhD thesis in the mid-1960s
(although the published version only 
appeared much later in~\cite{Verdier96}). Not surprisingly, the idea
of studying the finite objects in these categories followed suit 
soon after, see Illusie~\cite{Illusie71A,Illusie71B,Illusie71C}.

Right from the start there was a pervasive discomfort with derived
and triangulated
categories---the intuition that had been built up, in dealing with
concrete categories, mostly fails for triangulated categories.
In case the reader is wondering: in the previous sentence
the word ``concrete'' has a precise, technical meaning, and it
is an old theorem of Freyd~\cite{Freyd70,Freyd04} that triangulated
categories often aren't concrete.
Further testimony, to the strangeness of derived and triangulated categories, is
that it took two decades before the intuitive notion of 
finiteness, which dates back to Illusie's 
articles~\cite{Illusie71A,Illusie71B,Illusie71C},
 was given its
correct formal definition. The following may be found in
\cite[Definition~1.1]{Neeman92A}.

\dfn{D1.1}
Let $\ct$ be a triangulated category with coproducts. An object $C\in\ct$ is 
called \emph{compact} if $\Hom(C,-)$ commutes with coproducts. The full subcategory
of all compact objects will be denoted by $\ct^c$.
\edfn

\rmk{R1.3}
I have often been asked where the name ``compact''
came from. In the preprint
version of \cite{Neeman92A} these objects went by a different
name, but the (anonymous) referee
didn't like it. I was given a choice: I was allowed to baptize 
them either ``compact''
or ``small''. 

Who was I to argue with a referee?
\ermk

Once one has a good working definition of what the finite objects ought to be,
the next step is to give the right criterion which guarantees that the category
has ``enough'' of them. For triangulated categories the right definition
didn't come until \cite[Definition~1.7]{Neeman96}.

\dfn{D1.5}
Let $\ct$ be a triangulated category with coproducts. The category $\ct$ is
called \emph{compactly generated} if every nonzero object $X\in\ct$ admits a
nonzero map $C\la X$, with $C\in\ct$ a compact object.
\edfn

As the reader may have guessed from the name, compactly generated triangulated
categories are ones in which it is often possible to reduce general problems 
to questions about compact objects---which tend to be easier.

All of the above nowadays counts as ``classical'', meaning that it 
is two or more decades old and there
is already a substantial and diverse literature exploiting the ideas.
This article explores the recent developments that arose from trying to
understand how efficiently one can approximate arbitrary objects by
compact ones.
We first survey 
the results obtained to date. This review is
on the skimpy side, partly because there already are other,
more expansive 
published accounts in
the literature, but mostly because we want to leave ourselves
space to suggest possible directions
for future research. Thus the article can be thought of as having two
components: a bare-bone review of what has been achieved to date,
occupying  Sections~\ref{S77} to \ref{S5}, 
followed by Section~\ref{S6} which is comprised of suggestions of
avenues that might merit further development. 

Our review presents 
just enough detail so that the open questions, making  
up Section~\ref{S6}, 
can be formulated clearly and comprehensibly, and so that their
significance and potential applications can be illuminated.
This has the 
unfortunate side effect that we give short shrift to the many deep,
substantial contributions,
made by numerous mathematicians, which preceded and inspired the work
presented here.
The author apologizes in advance for this omission, which is the 
inescapable corollary 
of page limits.
The reader is referred to the other surveys of the subject, where
more care is taken to attribute the ideas correctly to their
originators, and give credit where credit is due. 

We permit ourselves to gloss over difficult technicalities, 
nonchalantly skating 
by nuances and subtleties, with only an occasional 
passing reference to the other surveys or to the
research papers for more detail.

The reader wishing to begin with examples and applications, to keep in mind
through the forthcoming abstraction, is encouraged  to first
look at the Introduction to \cite{Neeman17B}.

\section{Approximable triangulated categories---the
formal definition as a variant on Fourier series}
\label{S77}

It's now time to start our review, offering 
a glimpse of the recent progress that
was made by trying to measure how ``complicated'' an object is, in other
words how far it is from being compact. What follows is sufficiently new 
for there to be much room for 
improvement: the future will undoubtedly see cleaner, 
more elegant and more 
general formulations. What's presented here is the current crude state
of this emerging field.

\dis{D1.6}
This section is devoted to defining approximable triangulated categories, and
the definition is technical and  at first sight could appear artificial, 
maybe even forbidding. 
It might help
therefore to motivate it with an analogy.

Let $\sS^1$ be the circle, and let $M(\sS^1)$ be the set of all 
complex-valued,
Lebesgue-measurable functions on $\sS^1$. As
usual we view $\sS^1=\rr/\zz$ as
the quotient of its universal cover $\rr$ by the fundamental group $\zz$; 
this identifies functions
on $\sS^1$ with periodic functions on $\rr$ with period $1$. In 
particular the function $g(x)=e^{2\pi ix}$ belongs to $M(\sS^1)$.
And, for each $\ell\in\zz$, we have that $g(x)^\ell=e^{2\pi i\ell x}$ also
belongs to $M(\sS^1)$. Given a norm on 
the space $M(\sS^1)$, for example the $L^p$--norm,
we can try to approximate arbitrary $f\in M(\sS^1)$ by 
Laurent polynomials in $g$, that is look for complex numbers
$\{\lambda_\ell^{}\in\cC\mid-n\leq \ell\leq n\}$ such that
\[
\left\lVert f(x)-\sum_{\ell=-n}^n\lambda_\ell^{} g(x)^\ell\right\rVert_p
=
\left\lVert f(x)-\sum_{\ell=-n}^n\lambda_\ell^{} e^{2\pi i\ell x}\right\rVert_p
<\e
\]
with $\e>0$ small.
This leads us to the familiar territory of Fourier series.

Now imagine trying to do the same, but replacing $M(\sS^1)$ 
by a triangulated category.
Given a triangulated category $\ct$, which we assume
to have coproducts, we would like to pretend to do Fourier analysis 
on it. We would need to choose:
\begin{enumerate}
\item
Some analog of the function $g(x)=e^{2\pi ix}$. Our replacement
for this will be to choose a compact generator $G\in\ct$.
Recall: a \emph{compact generator} is a compact object
$G\in\ct$, such that every nonzero object $X\in\ct$ admits a
nonzero map $G[i]\la X$ for some $i\in\zz$.
\item
We need to choose something like a metric, the analog of the 
$L^p$--norm on $M(\sS^1)$. For us this will be done by
picking a t-structure $\left(\ct^{\leq0},\ct^{\geq0}\right)$ on $\ct$.
The heuristic is that we will view a morphism $E\la F$ in $\ct$ as
``short'' if, in the triangle $E\la F\la D$, the object $D$ belongs
to $\ct^{\leq-n}$ for large $n$. We will come
back to this in Discussion~\ref{D5.17}.
\item
We need to have an analog of the construction that passes,
from the function $g(x)=e^{2\pi ix}$ and the integer
$n>0$, to the vector space of trigonometric Laurent polynomials
$\sum_{\ell=-n}^n\lambda_\ell^{} e^{2\pi i\ell x}$.
\end{enumerate}
As it happens our solution to (3) is technical. We need a recipe, that
begins with the object $G$ and the integer $n>0$, and proceeds to cook up
a collection of  more objects. We
ask the reader
to accept it as a black box, with only a sketchy explanation 
just before Remark~\ref{R1.73}.
\edis

\bbx{BB1.7}
Let $\ct$ be a triangulated category and let $G\in\ct$ be an object.
Let $n>0$ be an integer.
We will have occasion to refer to the following
four full subcategories of $\ct$.
\begin{enumerate}
\item
The subcategory $\genuf Gn\subset\ct$ is defined unconditionally, 
and if $\ct$ has coproducts one can also define 
the larger subcategory $\ogenuf Gn$. Both of these
subcategories are classical, the reader can find 
the subcategory $\genuf Gn$ in
Bondal and Van den 
Bergh~\cite[the discussion between 
Lemma~2.2.2 and Definition~2.2.3]{BondalvandenBergh04},
and the subcategory $\ogenuf Gn$ 
in~\cite[the discussion between Definition~2.2.3
and Proposition~2.2.4]{BondalvandenBergh04}.
\item
If the category $\ct$ has coproducts, we will also have occasion
to consider the full subcategory $\ogenul G{}n$. Once again this
category is classical (although the name isn't). The reader can
find it in Alonso, Jerem{\'{\i}}as and Souto~\cite{Alonso-Jeremias-Souto03},
where it would go 
by the name ``the cocomplete pre-aisle generated by $G[-n]$''.
\item
Once again assume that $\ct$ has coproducts. Then we will also
look at the full subcategory $\ogenu Gn{-n,n}$. This 
construction is relatively new.
\end{enumerate} 
Below we give a vague description of what is going on in these
constructions; but when it comes to the technicalities we
ask the reader to either
accept these as black boxes, or
refer to \cite[Reminder~0.8~(vii), (xi) and (xii)]{Neeman17A}
for detail. We mention that there is a slight clash of notation
in the literature: 
what we call 
$\ogenuf Gn$ in (1), following Bondal and Van den Bergh, goes by
a different name in 
\cite[Reminder~0.8~(xi)]{Neeman17A}. The name it goes by there
is the case $A=-\infty$ and $B=\infty$ of the more general subcategory
$\ogenu Gn{A,B}$.

Now for the vague explanation of what goes on
in (1), (2) and (3) above. In a triangulated category $\ct$, there aren't many ways to build
new objects out of old ones. One can shift objects, form direct summands,
form finite
direct sums (or infinite ones if coproducts exist), and one can form
extensions. In the categories $\genuf Gn$
and $\ogenuf Gn$ of (1) there is a bound on the number
of allowed extensions, and the difference 
between the two is whether infinite coproducts are allowed. 
In the category $\ogenul G{}n$ of (2) the 
bound is on the 
permitted shifts. And in the category $\ogenu Gn{-n,n}$ of (3), both the shifts 
allowed and the number of extensions permitted are restricted.
\ebbx

\rmk{R1.73}
The reader should note that an example would not be illuminating,
the categories $\genuf Gn$, $\ogenuf Gn$, $\ogenul G{}n$ and $\ogenu Gn{-n,n}$
are not usually overly computable. For example: 
let $R$ be an associative ring, and let $\ct=\D(R)$ be the
unbounded derived category of complexes of left $R$-modules. The object
$R\in\ct$, that is the complex which is $R$ in degree zero and vanishes 
in all other degrees, is a compact generator for $\ct=\D(R)$.

But if we wonder what the categories $\genuf Rn$,
$\ogenuf Rn$, $\ogenul R{}n$ and 
$\ogenu Rn{-n,n}$ might turn out to be, only the category
$\ogenul R{}n$ is straightforward: it is the category of all
cochain complexes whose cohomology vanishes in degrees $>n$.
The three categories $\genuf Rn$, $\ogenuf Rn$ and $\ogenu Rn{-n,n}$
are mysterious in general. In fact, the computation of
$\genuf Gn$ is the subject of conjectures that have attracted much
interest. We will say a tiny bit about theorems in this direction
in Section~\ref{S3}, and will mention one of the 
active, open conjectures
in the discussion between Definition~\ref{D6.11} and Problem~\ref{P6.13}.
\ermk

\rmk{R1.95}
In the definition of approximable triangulated categories, which is about
to come, the category $\ogenu Gn{-n,n}$ will play the role of the 
replacement for the vector space of trigonometric Laurent polynomials of 
degree $\leq n$,
which came up in the desiderata of
Discussion~\ref{D1.6}(3). The older categories 
$\genuf Gn$, $\ogenuf Gn$ and $\ogenul G{}n$ will be needed 
later in the article.
\ermk

\rmk{R1.97}
Let us return to the heuristics of Discussion~\ref{D1.6}. Assume we have chosen
the t-structure $\left(\ct^{\leq0},\ct^{\geq0}\right)$ as in
Discussion~\ref{D1.6}(2), which we think of as our replacement for the
$L^p$--norm on $M(\sS^1)$. And we have also chosen a compact generator
$G\in\ct$ as in Discussion~\ref{D1.6}(1), which we think of as the
analog of the exponential function $g(x)=e^{2\pi ix}$. We have declared that
the subcategories $\ogenu Gn{-n,n}$ will be our replacement for
the vector space of trigonometric Laurent
polynomials of degree $\leq n$,
as in Discussion~\ref{D1.6}(3). It's now
time to start approximating functions by trigonometric Laurent 
polynomials. 

Let us therefore assume
we start with some object $F\in\ct$, and find a good approximation of
it by the object $E\in\ogenu Gm{-m,m}$. Recall: this means that we find
a morphism $E\la F$ such that, in the triangle $E\la F\la D$, the
object $D$ belongs to $\ct^{\leq-M}$ for some suitably large $M$.

Now we can try to iterate, and find a good approximation for $D$. Thus 
we can look for a morphism $E''\la D$, with $E''\in\ogenu Gn{-n,n}$,
and such that in the triangle $E''\la D\la D'$ the object $D'$ belongs to
$\ct^{\leq-N}$, with $N>M$ even more enormous than $M$. Can we combine
these to improve our initial approximation of $F$?

To do this, let's build up 
the octahedron on the composable morphisms $F\la D\la D'$. We end up
with a diagram where the rows and columns are triangles
\[\xymatrix{
E\ar[r]\ar@{=}[d] & E'\ar[r]\ar[d] & E''\ar[d] \\
E \ar[r] & F\ar[r]\ar[d] & D\ar[d] \\
      & D'\ar@{=}[r] & D'
}\]
and in particular the triangle $E'\la F\la D'$ gives that
$E'$ is an even better approximation of $F$ than $E$ was.
We are therefore interested in knowing if the triangle $E\la E'\la E''$,
coupled with the fact that $E\in\ogenu Gm{-m,m}$ and $E''\in\ogenu Gn{-n,n}$,
gives any information about where $E'$ might lie with
respect to the construction of Black Box~\ref{BB1.7}(3). Hence it is
useful to know the following.
\ermk

\fac{F1.8}
Let $\ct$ be a triangulated category with coproducts.
The construction of Black Box~\ref{BB1.7}(3) satisfies
\begin{enumerate}
\item
If $E$ is an object of $\ogenu Gn{-n,n}$, then the shifts $E[1]$ and
$E[-1]$ both belong to
$\ogenu G{n+1}{-n-1,n+1}$.
\item
Given an exact triangle 
$E\la E'\la E''$, with 
 $E\in\ogenu Gm{-m,m}$ and $E''\in\ogenu Gn{-n,n}$, it
follows that 
$E'\in\ogenu G{m+n}{-m-n,m+n}$.
\end{enumerate}
\efac

Combining Remark~\ref{R1.97} with Facts~\ref{F1.8} allows us to improve
approximations through iteration. 
Hence part (2) of the definition below 
becomes natural, it iterates to
provide arbitrarily good approximations.

\dfn{D1.11}
Let $\ct$ be a triangulated category with coproducts. It is \emph{approximable}
if there exist 
a t-structure $\left(\ct^{\leq0},\ct^{\geq0}\right)$,
a compact generator $G\in\ct$, and an integer $n>0$ such that
\begin{enumerate}
\item
$G$ belongs to $\ct^{\leq n}$ and $\Hom\left(G,\ct^{\leq-n}\right)=0$.
\item
Every object $X\in\ct^{\leq0}$ admits an exact triangle 
$E\la X\la D$ with $E\in\ogenu Gn{-n,n}$ and 
with $D\in\ct^{\leq-1}$.
\end{enumerate}
\edfn

\rmk{R1.15}
While part (2) of Definition~\ref{D1.11} comes motivated by the analogy
with Fourier analysis, part (1) of the definition seems random.
It requires the t-structure, which is our replacement for the $L^p$--norm,
to be compatible with the compact generator, which is the analog of
$g(x)=e^{2\pi i x}$. As the reader will see in Proposition~\ref{P4.9},
this has the effect of uniquely specifying the t-structure (up
to equivalence). So maybe a better parallel would be to fix our norm to 
be a particularly nice one, for example the $L^2$--norm on $M(\sS^1)$.

Let me repeat myself: as with all new mathematics, Definition~\ref{D1.11}
should be viewed as provisional. In
the remainder of this survey we will discuss the applications
as they now stand, to highlight the power of the methods. But
I wouldn't be surprised in the slightest
if future applications turn out to require modifications,
and/or generalizations, of the definitions and of the theorems
that have worked so far.
\ermk

\section{Examples of approximable 
triangulated categories}
\label{S2}

In Section~\ref{S1} we gave the definition of approximable triangulated
categories. The definition combines old, classical ingredients (t-structures
and compact generators) with a new construction, the category
$\ogenu Gn{-n,n}$ of Black Box~\ref{BB1.7}(3). The first thing to show
is that the theory is nonempty: we need to produce examples,
categories people care about which satisfy the definition of approximability.
The current section is devoted to the known examples of approximable
triangulated categories. We repeat what we have said before: the subject
is in its infancy, there could well be many more examples out there.

\exm{E2.1}
Let $\ct$ be a triangulated category with coproducts. If $G\in\ct$ is a
compact generator, such that $\Hom(G,G[i])=0$ for all $i>0$, then 
the category $\ct$ is approximable.

This example turns out to be easy, the reader is
referred to \cite[Example~3.3]{Neeman17A} for the (short) proof. 
Special cases include 
\begin{enumerate}
\item
$\ct=\D(\Mod R)$, where $R$ is a dga with $H^i(R)=0$ for $i>0$.
\item
The homotopy category of spectra.
\end{enumerate}
\eexm

\exm{E2.3}
If $X$ is a quasicompact, separated scheme, then the category
$\Dqc(X)$ is approximable. We remind the reader of the traditional
notation being used here: the category 
$\D(X)$ is the unbounded derived
category of complexes of sheaves of $\co_X^{}$--modules, and
the full subcategory $\Dqc(X)\subset\D(X)$ has for objects the
complexes with quasicoherent cohomology.

The proof of the approximability of $\Dqc(X)$ isn't trivial. The
category has a standard t-structure, that part is easy. The 
existence of a compact generator $G$ needs proof, it may be
found in Bondal and Van den 
Bergh~\cite[Theorem~3.1.1(ii)]{BondalvandenBergh04}. Their proof 
isn't constructive, it's only an existence proof, but 
it does give enough
information to deduce that 
part (1) of Definition~\ref{D1.11} is satisfied
by every compact generator (indeed, it's satisfied by every
compact object). See
\cite[Theorem~3.1.1(i)]{BondalvandenBergh04}.
But 
it is a challenge to show that we may choose a compact
generator $G$ and an integer $n>0$ in such a way that 
Definition~\ref{D1.11}(2) is satisfied.

If we further 
assume that $X$ is of finite type over a noetherian ring $R$,
then the (relatively intricate) 
proof of the approximability of $\Dqc(X)$
occupies
\cite[Sections~4 and 5]{Neeman17}. The little trick, that extends 
the result to all quasicompact and separated $X$, wasn't
observed until later: it appears  
in \cite[Lemma~3.5]{Neeman17A}.
\eexm

\exm{E2.5}
It is a theorem that, under mild hypotheses, the recollement of any
two approximable triangulated categories is approximable.
To state the ``mild hypotheses'' precisely:
suppose we are given a recollement of triangulated categories
\[\xymatrix@C+20pt{
\car\ar[r] &
\cs\ar@<0.7ex>[l]\ar@<-0.7ex>[l]\ar[r]  &
\ct\ar@<0.7ex>[l]\ar@<-0.7ex>[l]
}\]
with $\car$ and $\ct$ approximable. Assume further that
the category $\cs$ is compactly generated, and any compact object
$H\in\cs$ has the property that $\Hom\big(H,H[i]\big)=0$ for $i\gg0$.
Then the category $\cs$ is also approximable.

The reader can find
the proof in~\cite[Theorem~4.1]{Burke-Neeman-Pauwels18},
it is the main result in the paper.
The bulk of the article is devoted to developing the machinery
necessary to prove the theorem---hence it's worth noting 
that this machinery
has since demonstrated usefulness in other contexts, see the
subsequent
articles~\cite{Neeman18,Neeman18A}.

There is a beautiful theory of noncommutative schemes, and a rich literature
studying them. And many of the interesting examples of such schemes are obtained
as recollements of ordinary schemes, or of admissible pieces of them. Thus
the theorem that recollements of approximable triangulated categories
are approximable gives a wealth of new examples of approximable triangulated
categories. 

Since this ICM is being held in St Petersburg, it 
would be remiss not to mention that the theory of noncommutative algebraic
geometry, in the sense of the previous paragraph, is a subject to which
Russian mathematicians have contributed a vast amount. The 
seminal work
of Bondal, Kontsevich, Kuznetsov, Lunts and Orlov immediately springs
to mind. For a beautiful introduction to the field the reader might
wish to look at the early sections of Orlov~\cite{Orlov16}. The later
sections prove an amazing new theorem, but the early ones give a lovely
survey of the background. In fact: the theory sketched in this survey was
born when I was trying to read and understand Orlov's beautiful article.
\eexm

\section{Applications: strong generation}
\label{S3} 

We begin by reminding the reader of a classical definition, going back
to Bondal and Van den Bergh~\cite{BondalvandenBergh04}.

\dfn{D3.1} Let $\ct$ be
triangulated category. An object $G\in\ct$ is 
called a \emph{strong generator} if  there exists an integer
$\ell>0$ with $\ct=\genuf G\ell$, where the notation is as
in Black Box~\ref{BB1.7}(1). The
category $\ct$ is called \emph{regular} or \emph{strongly generated}
if it contains a strong generator.
\edfn 

The first application of approximability is the proof of the following
two theorems.

\thm{T3.3}
Let $X$ be a quasicompact, separated scheme. 
The derived category of perfect complexes
on $X$, denoted here by $\dperf X$, is regular if and only if $X$ has
a cover by open subsets $\spec{R_i}\subset X$, with each $R_i$ of finite
global dimension.
\ethm

\rmk{R3.5}
If $X$ is noetherian and separated, then Theorem~\ref{T3.3}
specializes to saying that $\dperf X$ is regular if and
only if $X$ is regular and finite dimensional. Hence the terminology.
\ermk

\thm{T3.7}
Let $X$ be a noetherian, separated, finite-dimensional, 
quasiexcellent scheme. Then the category
$\D^b\big(\mathrm{Coh}(X)\big)$, the bounded derived category of coherent 
sheaves on $X$, is always regular. 
\ethm

\rmk{R3.9}
The reader is referred to \cite{Neeman17} and to Aoki~\cite{Aoki20} 
for the proofs of Theorems~\ref{T3.3} and \ref{T3.7}. More precisely:
for Theorem~\ref{T3.3} see \cite[Theorem~0.5]{Neeman17}. About
Theorem~\ref{T3.7}: if we add the assumption that every closed
subvariety of $X$ admits a regular alteration then the result
may be found in \cite[Theorem~0.15]{Neeman17}, but Aoki~\cite{Aoki20}
found a lovely argument that allowed him to extend the statement
to all quasiexcellent $X$.

There is a rich literature on strong generation, with beautiful
papers by many authors. In the introduction to \cite{Neeman17},
as well as in \cite{Neeman16} and
\cite[Section~7]{Neeman17B}, the reader can find an extensive
discussion of (some of) this fascinating work and of the way Theorems
\ref{T3.3} and \ref{T3.7} compare to the older literature.
For a survey taking an entirely different tack
see Minami~\cite{Minami21}, which places in historical perspective a 
couple of the key steps in the proofs that~\cite{Neeman17}
gives for Theorems
\ref{T3.3} and \ref{T3.7}.

Since all of this is now well documented in the published literature,
let us focus the remainder of the current survey
on the other applications of approximability. Those
are all 
still in preprint form, see~\cite{Neeman17A,Neeman18,Neeman18A},
although there are 
(published) surveys in \cite[Sections~8 and 9]{Neeman17B} and
in~\cite{Neeman19}. 
Those surveys
are fuller and more complete than the sketchy  one we are about to
embark on. As we present the material, we will feel free to refer the reader 
to the more extensive surveys whenever 
we deem it appropriate.
\ermk

\section{The freedom in the choice of compact generator
and t-structure}
\label{S4}

Definition~\ref{D1.11} tells us that a triangulated category $\ct$
with coproducts is
approximable if there exist, in $\ct$, a compact generator $G$ and
a t-structure $\left(\ct^{\leq0},\ct^{\geq0}\right)$ satisfying some
properties. The time has come to explore just how free we are
in the choice of the compact generator and of the t-structure. 
To address this question we begin by formulating:

\dfn{D4.5}
Let $\ct$ be a triangulated category. Two t-structures 
$\left(\ct_1^{\leq0},\ct_1^{\geq0}\right)$ and 
$\left(\ct_2^{\leq0},\ct_2^{\geq0}\right)$ are declared \emph{equivalent}
if there exists an integer $n>0$ such that
\[
\ct_1^{\leq-n}\subset\ct_2^{\leq0}\subset\ct_1^{\leq n}.
\]
\edfn

\dis{D4.7}
Let $\ct$ be a triangulated category with coproducts. If $G\in\ct$ is a
compact object and $\ogenul G{}0$ is as in Black Box~\ref{BB1.7}(2),
then Alonso, Jerem{\'{\i}}as and 
Souto~\cite[Theorem~A.1]{Alonso-Jeremias-Souto03}, building
on the work of Keller and Vossieck~\cite{Keller-Vossieck88}, teaches us that
there is a unique t-structure $\left(\ct^{\leq0},\ct^{\geq0}\right)$
with $\ct^{\leq 0}=\ogenul G{}n$.
We will call this the \emph{t-structure
generated by $G$,} and denote it $\left(\ct_G^{\leq0},\ct_G^{\geq0}\right)$. 

In Black Box~\ref{BB1.7}(2) we asked the reader to 
accept, as a black box,
the construction passing from an
object $G\in\ct$ to the
subcategory $\ogenul G{}{0}$. If $G$ is compact, 
then~\cite[Theorem~A.1]{Alonso-Jeremias-Souto03}
allows us to express this as $\ct^{\leq0}_G$ for a unique
t-structure. We ask the reader to 
accept on faith
that:
\edis

\lem{L4.89}
If $G$ and $H$ are two compact generators for the triangulated 
category $\ct$, then the two t-structures 
$\left(\ct_G^{\leq0},\ct_G^{\geq0}\right)$ and 
$\left(\ct_H^{\leq0},\ct_H^{\geq0}\right)$ are equivalent
as in Definition~\ref{D4.7}.
\elem

As it happens the proof of Lemma~\ref{L4.89} is easy,
the interested reader can find it in \cite[Remark~0.15]{Neeman17A}.
And Lemma~\ref{L4.89} leads us to:

\dfn{D4.95}
Let 
$\ct$ be a triangulated category in which there
exists a compact generator.  We define the \emph{preferred equivalence
class of t-structures} as follows: a t-structure belongs to the 
preferred equivalence
class if it is equivalent to $\left(\ct_G^{\leq0},\ct_G^{\geq0}\right)$ for some
compact generator $G\in\ct$, and by
Lemma~\ref{L4.89} it is equivalent to
 $\left(\ct_H^{\leq0},\ct_H^{\geq0}\right)$ for every compact generator
$H$.
\edfn

The following is also not too hard, and
may be found in \cite[Propositions~2.4 and 2.6]{Neeman17A}.

\pro{P4.9}
Let $\ct$ be an approximable triangulated category. Then for any t-stucture
$\left(\ct^{\leq0},\ct^{\geq0}\right)$ in the preferred
equivalence class, and
for any compact generator $H\in\ct$, there exists an integer
$n>0$ (which may depend on $H$ and on the t-structure), satisfying
\begin{enumerate}
\item
$H$ belongs to $\ct^{\leq n}$ and $\Hom\left(H,\ct^{\leq-n}\right)=0$.
\item
Every object $X\in\ct^{\leq0}$ admits an exact triangle 
$E\la X\la D$ with $E\in\ogenu Hn{-n,n}$ and
 with $D\in\ct^{\leq-1}$.
\end{enumerate}
Moreover: if $H$ is a compact generator, if $\left(\ct^{\leq0},\ct^{\geq0}\right)$
is a t-structure, and if there exists an integer $n>0$ satisfying (1) and (2)
above, then the t-structure 
$\left(\ct^{\leq0},\ct^{\geq0}\right)$
must belong to the preferred equivalence class.
\epro

\rmk{R4.11}
Strangely enough, the value of Proposition~\ref{P4.9} can be that it allows us
to find an explicit t-structure in the preferred equivalence class. 

Consider the case where $X$ is a quasicompact, separated scheme. By
Bondal and Van den 
Bergh~\cite[Theorem~3.1.1(ii)]{BondalvandenBergh04} we know that the category
$\Dqc(X)$ has a compact generator, but in Example~\ref{E2.3} we mentioned
that the existence proof isn't overly constructive, it doesn't give us
a handle on any explicit compact generator. Let $G$ be some
compact generator. From
Alonso, Jerem{\'{\i}}as and Souto~\cite[Theorem~A.1]{Alonso-Jeremias-Souto03}
we know that 
the subcategory $\ogenul G{}0$ of Black Box~\ref{BB1.7}(2) is equal
to $\ct_G^{\leq0}$ for a unique t-structure
$\left(\ct_G^{\leq0},\ct_G^{\geq0}\right)$ in the preferred equivalence class. 
But this doesn't leave us 
a whole lot wiser---the compact generator $G$ isn't explicit, hence
neither is the t-structure.

However: the combination of \cite[Theorem~5.8]{Neeman17} and 
\cite[Lemma~3.5]{Neeman17A} tells us that the category $\Dqc(X)$ is approximable,
and it so happens that the t-structure \emph{used in the proof,}
that is the t-structure for which a compact generator
$H$ and an integer $n>0$ satisfying (1) and (2) of Proposition~\ref{P4.9} are
shown to exist, happens to be the standard t-structure. From 
Proposition~\ref{P4.9} we now deduce that the standard t-structure is in
the preferred equivalence class.
\ermk

\section{Structure theorems in approximable triangulated
categories}
\label{S5}

An approximable triangulated category $\ct$ must have a compact generator $G$, 
and
Definition~\ref{D4.95} constructed for us a preferred equivalence class of 
t-structures---namely
all those equivalent to  $\left(\ct_G^{\leq0},\ct_G^{\geq0}\right)$. 
Recall that,
for any t-structure $\left(\ct^{\leq0},\ct^{\geq0}\right)$, it is
customary to define
\[
\ct^-=\cup_{n=1}^\infty\ct^{\leq n},\qquad
\ct^+=\cup_{n=1}^\infty\ct^{\geq-n},\qquad
\ct^b=\ct^-\cap\ct^+.
\]
It is an easy exercise to show, directly from Definition~\ref{D4.5}, that equivalent
t-structures give rise to identical $\ct^-$, $\ct^+$ and $\ct^b$. Therefore
triangulated categories with a single compact generator, and in particular
approximable triangulated categories, have preferred subcategories
$\ct^-$, $\ct^+$ and $\ct^b$, which are intrinsic---they are simply the ones 
corresponding to any t-structure in the preferred
equivalence class. In the remainder of this survey we will assume 
that $\ct^-$, $\ct^+$ and $\ct^b$ always stand for the preferred ones.

In the heuristics of Discussion~\ref{D1.6}(2) we told the reader that 
a t-structure $\left(\ct^{\leq0},\ct^{\geq0}\right)$ is to be viewed as
a metric on $\ct$. In Definition~\ref{D5.3} below, the heuristic
is that we construct a full subcategory $\ct^-_c$ to be the closure
of $\ct^c$ with respect to any of the (equivalent) metrics that come
from t-structures in the preferred equivalence class.

\dfn{D5.3}
Let $\ct$ be an approximable triangulated category. The full
subcategory $\ct^-_c$ is given by
\[
\text{Ob}(\ct^-_c)=
\left\{F\in\ct\,\,\left|\,\,
\begin{array}{c}
\text{For every integer $n>0$ and for every t-structure}\\
\left(\ct^{\leq0},\ct^{\geq0}\right)\text{ in the preferred equivalence class,}\\
\text{there exists an exact triangle $E\la F\la D$ in $\ct$}\\
\text{with $E\in\ct^c$ and $D\in\ct^{\leq-n}$} 
\end{array}
\right.
\right\}.
\]
The full subcategory $\ct^b_c$ is defined to be $\ct^b_c=\ct^-_c\cap\ct^b$.
\edfn

\rmk{R5.99999}
Let $\ct$ be an approximable triangulated category. Aside from the classical,
full
subcategory $\ct^c$ of compact objects, which we encountered back in
Definition~\ref{D1.1}, we have in this section concocted five more 
intrinsic,
full  subcategories of $\ct$: they are 
$\ct^-$, $\ct^+$, $\ct^b$, $\ct^-_c$ and $\ct^b_c$. It can be proved
that all six subcategories, that is the old $\ct^c$ and the five new
ones, are thick subcategories of $\ct$. In particular each of them
is a triangulated category.
\ermk

\exm{E5.7}
It becomes
interesting to figure out what all these categories
come down to in examples.

Let $X$ be a quasicompact, separated scheme.
From Example~\ref{E2.3} we know that the category $\ct=\Dqc(X)$ 
is approximable, and in Remark~\ref{R4.11} we noted that the
standard t-structure is in the preferred equivalence class. 
This can be used to show that, for $\ct=\Dqc(X)$, we have
\[
\begin{array}{rclcrclcrcl}
\ct^-&=&\Dqcmi(X),&\,\qquad\,&\ct^+&=&\Dqcpl(X),&\,\qquad\,&\ct^b&=&\Dqcb(X),\\
\ct^c&=&\dperf X,&\,\qquad\,&\ct^-_c&=&\dcohm(X),&\,\qquad\,&\ct^b_c&=&\dcoh(X),
\end{array}
\]
where the last two equalities
assume that the scheme $X$ is noetherian,
and all six categories on the right of the equalities 
have their traditional meanings.

The reader can find an extensive discussion of the claims above in
\cite{Neeman17B}, more precisely in the paragraphs between
\cite[Proposition~8.10]{Neeman17B} and \cite[Theorem~8.16]{Neeman17B}.
That discussion goes beyond the scope of the current survey,
it analyzes the categories $\ct^b_c\subset\ct^-_c$ in the generality
of
non-noetherian schemes, where they still have a classical 
description---of course not involving the category of coherent sheaves.
After all coherent sheaves 
do not behave well for non-noetherian schemes.
\eexm

\rmk{R5.9}
In this survey we spent some effort introducing the notion of
approximable triangulated categories. In Example~\ref{E2.3}
we told the reader that it is a theorem (and not a trivial one)
that, as long as a scheme $X$ is quasicompact and separated,
the derived category $\Dqc(X)$ is approximable. In this section
we showed that every approximable triangulated category comes with
canonically defined, intrinsic subcategories 
$\ct^-$, $\ct^+$, $\ct^b$, $\ct^c$, $\ct^-_c$ and $\ct^b_c$, and
in Example~\ref{E5.7} we informed the reader that, in the
special case where $\ct=\Dqc(X)$, these turn out
to be, respectively, $\Dqcmi(X)$, $\Dqcpl(X)$, $\Dqcb(X)$,
$\dperf X$, $\dcohm(X)$ and $\dcoh(X)$.

Big deal. This teaches us that the traditional subcategories
$\Dqcmi(X)$, $\Dqcpl(X)$, $\Dqcb(X)$,
$\dperf X$, $\dcohm(X)$ and $\dcoh(X)$ of the category $\Dqc(X)$
all have intrinsic descriptions. This might pass as a curiosity,
unless we can actually use it to prove something
we care about that we didn't use
to know.
\ermk

\dis{D5.909}
To motivate the next theorem, it might help to think of the parallel
with functional analysis.

Let $M(\rr)$ be the vector space of Lebesgue-measurable, real-valued functions
on $\rr$. Given any two functions $f,g\in M(\rr)$ we can pair them
by integrating the product, that is we form the pairing
\[
\langle f,g\rangle=\int fg\,d\mu
\]
where $\mu$ is Lebesgue measure. This gives us a map
\[\xymatrix@C+50pt{
M(\rr)\times M(\rr)
\ar[r]^-{\langle -,-\rangle} & \rr\cup\{\infty\},
}\]
where the integral is declared to be infinite if it doesn't converge.

We can restrict this pairing to subspaces of $M(\rr)$. For example
if $f\in L^p(\rr)$ and $g\in L^q(\rr)$ with $\frac1p+\frac1q=1$ then
the integral converges, that is $\langle f,g\rangle\in\rr$,
and we deduce a map
\[\xymatrix@C+20pt{
L^p(\rr)
\ar[r] & \Hom\big(L^q(\rr),\rr\big)
}\]
which turns out to be an isometry of Banach spaces. 

The category-theoretic version is that on any category $\ct$ there is 
the pairing sending two objects $A,B\in\ct$ to $\Hom(A,B)$. Of course
this pairing isn't symmetric, we have to keep track of the position
of $A$ and of $B$ in $\Hom(A,B)$. If $R$
is a commutative ring and $\ct$ happens to be an $R$-linear category,
then $\Hom(A,B)$ is an $R$-module and the pairing delivers a map
\[\xymatrix@C+50pt{
\ct\op\times \ct
\ar[r]^-{\Hom( -,-)} & \Mod R,
}\]
where the op keeps track of the variable in the first position.
And now we are free to restrict to subcategories of $\ct$.

If $\ct$ happens to be approximable and $R$-linear, 
we have just learned that it comes
with six intrinsic subcategories
$\ct^-$, $\ct^+$, $\ct^b$, $\ct^c$, $\ct^-_c$ and $\ct^b_c$.
We are free to restrict the Hom pairing to any couple of them. This gives 
us 36 possible pairings, and each of those yields two maps from a subcategory
to the dual of another. There are 72 cases we could study, and the 
theorem below tells us something useful about four of those.
\edis

\thm{T29.27}
Let $R$ be a noetherian ring, and let
$\ct$ be an $R$--linear, approximable triangulated category.
Suppose there exists in $\ct$ a compact generator $G$ so that
$\Hom\big(G,G[n]\big)$ is a finite $R$--module for all $n\in\zz$.
Consider the two functors
\[
\cy:\ct^-_c\la\Hom_R^{}\big([\ct^c]\op\,,\,\Mod R\big),\qquad
\wt\cy:\big[\ct^-_c\big]\op\la\Hom_R^{}\big(\ct^b_c\,,\,\Mod R\big)
\]
defined by the formulas $\cy(B)=\Hom(-,B)$ and
$\wt\cy(A)=\Hom(A,-)$, as in Discussion~\ref{D5.909}.
Now consider the following composites
\[\xymatrix@C+20pt@R-20pt{
\ct^b_c \ar@{^{(}->}[r]^i & \ct^-_c
\ar[r]^-{\cy} &
\Hom_R^{}\big([\ct^c]\op\,,\,\Mod R\big) \\
\big[\ct^c\big]\op \ar@{^{(}->}[r]^{\wi} & \big[\ct^-_c\big]\op
\ar[r]^-{\wt\cy} &
\Hom_R^{}\big(\ct^b_c\,,\,\Mod R\big)
}\]
We assert:
\be
\item
  The functor $\cy$ is full, and the essential image consists
  of the
  locally finite homological functors [see Explanation~\ref{X29.102345}
    for the definition of locally finite functors]. The composite $\cy\circ i$
  is fully faithful, and the essential image consists of the
  finite homological functors [again: see Explanation~\ref{X29.102345} 
  for the definition].
\item
  With the notation as in Black Box~\ref{BB1.7}(1), 
assume\footnote{\label{foot2}
What's important for the current survey is that,
if $X$ is a noetherian, separated scheme, then
$\ct=\Dqc(X)$ satisfies this hypothesis 
provided $X$ is finite-dimensional and quasiexcellent.} 
that 
$\ct=\ogenuf {H}n$ for some integer $n>0$ and some object $H\in\ct^b_c$.
  Then
  the functor $\wt\cy$ is full, and the essential image consists
  of the
  locally finite homological functors. The composite $\wt\cy\circ \wi$
  is fully faithful, and the essential image consists of the
  finite homological functors.
\ee
\ethm

\xpl{X29.102345}
In the statement of Theorem~\ref{T29.27}, the \emph{locally finite}
functors, either of the form $H:\big[\ct^c\big]\op\la\Mod R$ or 
of the form $H:\ct^b_c\la\Mod R$, are the functors such that
\be
\item
$H\big(A[i]\big)$ is a finite $R$--module for every $i\in\zz$ and every $A$
in either $\ct^c$ or $\ct^b_c$.
\item
For fixed $A$, in one of $\ct^c$ or $\ct^b_c$, 
we have $H\big(A[i]\big)=0$ if $i\ll0$.
\setcounter{enumiv}{\value{enumi}}
\ee
The \emph{finite} functors are those for which we also have
\be
\setcounter{enumi}{\value{enumiv}}
\item
$H\big(A[i]\big)=0$ for all $i\gg0$.
\setcounter{enumiv}{\value{enumi}}
\ee
\expl

\rmk{R5.15}
The proof of part (1) of Theorem~\ref{T29.27} may be found in 
\cite{Neeman17A}, while the proof of part (2) 
of Theorem~\ref{T29.27} occupies \cite{Neeman18}. These
aren't easy theorems.

Let $\ct=\Dqc(X)$, with $X$ a scheme
proper over a noetherian ring $R$. Then the hypotheses of
Theorem~\ref{T29.27}(1) are satisfied. We learn 
(among other things) that
the natural functor, taking an object $B\in\dcoh(X)$ to the
$R$-linear functor
$\Hom(-,B):{\dperf X}\op\la\MMod R$, is a fully
faithful embedding
\[\xymatrix@C+50pt{
\dcoh(X)
\ar[r]^-{\cy\circ i} &
\Hom_R^{}\big({\dperf X}\op\,,\,\Mod R\big) 
}\] 
whose essential image is precisely the finite homological functors.

If we further assume that the scheme $X$ is finite-dimensional
and quasiexcellent then
the hypotheses of
Theorem~\ref{T29.27}(2) are also satisfied. We learn 
that the functor, taking an object $A\in\dperf X$ to the
$R$-linear functor $\Hom(A,-)$, is a fully faithful
embedding
\[\xymatrix@C+50pt{
{\dperf X}\op 
\ar[r]^-{\wt\cy\circ\wi} &
\Hom_R^{}\big(\dcoh(X)\,,\,\Mod R\big)
}\]
whose essential image is also the finite homological functors.

In \cite[Historical Survey~8.2]{Neeman17B} the reader
can find a discussion of the (algebro-geometric) precursors of
Theorem~\ref{T29.27}. As for the applications: let us go through
one of them.
\ermk

\rmk{R5.999}
Let $X$ be a scheme proper over the field $\cC$ of complex numbers,
and let $X\an$ be the underlying complex analytic space. The analytification
induces a functor we will call $\cl:\dcoh(X)\la\dcoh(X\an)$, it is the functor
taking a bounded complex of coherent algebraic sheaves on $X$ to the 
analytification, which is a bounded complex of coherent analytic 
sheaves on $X\an$. The pairing sending an object $A\in\dperf X$ 
and an object $B\in\dcoh(X\an)$ to $\Hom\big(\cl(A),B\big)$
delivers a map
\[\xymatrix@C+50pt{
\dcoh(X\an)
\ar[r] & \Hom_R^{}\big({\dperf X}\op\,,\,\Mod \cC\big).
}\]
Since the image lands in the finite homological functors, 
Theorem~\ref{T29.27}(1)
allows us to factor this uniquely through the inclusion $\cy\circ i$,
that is there exists 
(up to canonical natural isomorphism) a unique functor $\car$ 
rendering commutative the triangle
\[\xymatrix@C+50pt@R-20pt{
\dcoh(X\an)\ar@{.>}[dd]_{\exists!\car}
\ar[dr] & \\
 &\Hom_R^{}\big({\dperf X}\op\,,\,\Mod \cC\big).\\
\dcoh(X)\ar[ur]_-{\cy\circ i} &
}\]
And proving Serre's GAGA theorem reduces to the easy exercise of showing
that $\cl$ and $\car$ are inverse equivalences,
the reader can find this in
the (short) \cite[Section~8 and Appendix~A]{Neeman17A}.

The brilliant inspiration underpinning the approach is 
due to Jack Hall~\cite{Hall18}, 
he is the
person who came up with the idea of using the pairing
above, coupled with representability theorems, to prove GAGA.
The representability theorems available to Jack Hall at
the time weren't powerful
enough, and Theorem~\ref{T29.27} was motivated by trying to find
a direct path from the ingenious, simple idea to a fullblown proof.
\ermk

\dis{D5.17}
In preparation for the next theorem we give a very brief review
of metrics in triangulated categories. The reader is referred to
the survey article~\cite{Neeman19} for a much fuller and
more thorough account.

Given a triangulated category $\ct$, a \emph{metric} on $\ct$ assigns
a length to every morphism. In this article
the only metrics we consider are the ones arising from 
t-structures. If $\ct$ is an approximable triangulated category 
we choose a t-structure $\left(\ct^{\leq0},\ct^{\geq0}\right)$ 
in the preferred equivalence class, and this induces a metric as
follows.
Given a morphism $f:X\la Y$ we may complete to an exact triangle
$X\stackrel f\la Y\la D$, and the 
length of $f$ is given by the formula
\[
\text{Length}(f)=
\inf\left\{\left.\frac1{2^n}\,\,\right|\,
 n\in\zz\text{ and }D\in\ct^{\leq-n}\right\}.
\]
In this survey we allow the length of a morphism to be 
infinite; if the set on the right is empty then
we declare $\text{Length}(f)=\infty$.

This metric depends on the choice of t-structure, but not
much. As all
t-structures in the preferred equivalence class are equivalent,
any two preferred t-structures will give rise to equivalent metrics (with
an obvious definition of equivalence of metrics).

Note that if $\ct$ is a triangulated category and $\cs$ is a triangulated
subcategory, then a metric on $\ct$ restricts to a metric on $\cs$.
In particular: if 
$\ct$ is approximable, the metric on $\ct$ of the previous paragraph
restricts to give metrics on the full subcategories
$\ct^c$ and $\ct^b_c$. Once again these metrics are only defined
up to equivalence. And of course a metric on $\cs$ is also a metric
on $\cs\op$, thus we have specified (up to equivalence) canonical
metrics on
$\ct^c$, $\ct^b_c$, $\left[\ct^c\right]\op$ and  $\left[\ct^b_c\right]\op$.

Suppose $\cs$ is a triangulated category with a metric. 
A \emph{Cauchy sequence} in $\cs$ 
is a sequence of morphisms $E_1\la E_2\la E_3\la\cdots$ which
eventually become arbitrarily short. If $\ab$ is the category
of abelian groups, then the Yoneda embedding
$Y:\cs\la\MMod\cs$ embeds $\cs$ into the category $\MMod\cs$
of additive
functors $\cs\op\la\ab$. In the category $\MMod\cs$ colimits
exist, allowing us to define
\be
\item
The category $\fl(\cs)$ is the full subcategory of $\MMod\cs$, whose
objects are the colimits of Yoneda images
of Cauchy sequences in $\cs$.
\item
The full subcategory $\fs(\cs)\subset\fl(\cs)$ has for objects those
functors $F\in\fl(\cs)\subset\MMod\cs$ which take sufficiently short 
morphisms to isomorphisms. In symbols: $F\in\fl(\cs)$ belongs to
$\fs(\cs)$ if there exists an $\e>0$ such that
\[
\{\text{Length}(f)<\e\}\Longrightarrow\{F(f)\text{ is an isomorphism}\}.
\]
\item
The exact triangles in $\fs(\cs)$ are the colimits in $\MMod\cs$ of
Yoneda images of
Cauchy sequences of exact triangles in $\cs$,  
where the colimits happen to lie in 
$\fs(\cs)$. 
\ee
A word of
caution about (3): if we are given in $\cs$ 
a Cauchy sequence of exact triangles,
we can form the colimit in $\MMod\cs$ of its Yoneda image.
This colimit is guaranteed to lie in $\fl(\cs)$, but will
not usually lie in the smaller $\fs(\cs)$. If it happens to lie 
in $\fs(\cs)$ then (3) declares it to be an exact triangle in $\fs(\cs)$.
\edis

And now we are ready for the theorem.

\thm{T5.19}
Let $\cs$ be a triangulated category with a metric. Assume the metric is
\emph{good;} this is a technical term, see 
\cite[Definition~10]{Neeman19} for the precise formulation. Then
\be
\item
The category $\fs(\cs)$ of Discussion~\ref{D5.17}(2), with the exact
triangles as defined in Discussion~\ref{D5.17}(3), is a triangulated
category.
\setcounter{enumiv}{\value{enumi}}
\ee
Now let $\ct$ be an approximable triangulated category.
In Discussion~\ref{D5.17} we constructed (up to equivalence)
a metric on $\ct$, and hence on its subcategories $\ct^c$
and $\left[\ct^b_c\right]\op$. Those metrics are all good, 
and the theorem goes on to give natural, exact equivalences
of triangulated categories
\be
\setcounter{enumi}{\value{enumiv}}
\item
$\fs\left(\ct^c\right)\cong\ct^b_c$. This equivalence is unconditional.
\item
If the approximable triangulated category $\ct$ 
happens to be \emph{noetherian}
as in \cite[Definition~5.1]{Neeman18A}, then
$\fs\left(\left[\ct^b_c\right]\op\right)\cong\left[\ct^c\right]\op$.
\ee
\ethm

\rmk{R5.21}
First of all: in Theorem~\ref{T5.19}(3) we assumed that the approximable
triangulated $\ct$ is noetherian as in \cite[Definition~5.1]{Neeman18A}.
The only observation we want to make here is that if $X$ is a noetherian,
separated scheme then the approximable triangulated category
$\ct=\Dqc(X)$ is noetherian. Thus, for noetherian, separated schemes $X$,
Theorem~\ref{T5.19} gives exact equivalences of triangulated categories
\[
\fs\left(\dperf X\right)\cong\dcoh(X),\qquad
\fs\left(\dcoh(X)\op\right)\cong{\dperf X}\op.
\]
The research paper \cite{Neeman18A} contains the proofs of the
assertions 
in Theorem~\ref{T5.19}. The reader can
find a
skimpy survey in~\cite[Section~9]{Neeman17B}
and a more extensive one in \cite{Neeman19}. In 
\cite[Historical Survey~9.1]{Neeman17B} there is a discussion
of precursors of the results.
\ermk

\section{Future directions}
\label{S6}

New scientific developments are tentative and unpolished; 
only with the passage of time do they acquire the gloss and elegance
of a refined, varnished theory.
And there is nothing more difficult to predict than the future.
My colleague Neil Trudinger used to joke that my beard 
makes me look like a biblical
prophet---the reader should not be deceived, appearances
are notoriously misleading, the abundance of
facial hair isn't a reliable yardstick for measuring
the gift of foresight that marks out a visionary,
and I am certifiably not a 
clairvoyant. All I do in this section is offer a handful
of obvious
questions that spring to mind. The list is not meant to be
exhaustive, and might well be missing major 
tableaux of the overall picture. It
is entirely possible that the future will see this theory 
flourish in directions orthogonal to the ones
sketched here.

Let us begin with what's freshest in our minds: we have just seen
Theorem~\ref{T5.19}, part~(1) of
which tells us that, given a triangulated
category $\cs$ with a good metric, there is a recipe
producing another
triangulated category $\fs(\cs)$, which as it happens comes with
an induced
good metric. We can ask:

\plm{P6.1}
Can one formulate reasonable sufficient conditions, on the triangulated category
$\cs$ and on its good metric, to guarantee that
$\fs\left(\fs(\cs)\op\right)=\cs\op$? Who knows, maybe even necessary
and sufficient conditions?
\eplm

\mxm{M6.1.1}
Let
$\ct$ be an approximable triangulated
category and put $\cs=\ct^c$, with
the metric of Discussion~\ref{D5.17}. 
Theorem~\ref{T5.19}(2) computes for us that 
$\fs(\ct^c)\cong\ct^b_c$. I ask the reader to believe that the
natural, induced metric on $\fs(\ct^c)$ agrees with the metric 
on $\ct^b_c\subset\ct$ given in Discussion~\ref{D5.17}.
Now Theorem~\ref{T5.19}(3) goes on to tell us that, as long
as the approximable triangulated category $\ct$ is noetherian,
we also have that 
$\fs\left(\big[\ct^b_c\big]\op\right)\cong\big[\ct^c\big]\op$;
as it happens the induced good metric on
$\fs\left(\big[\ct^b_c\big]\op\right)$ also agrees, up to equivalence,
with the metric that Discussion~\ref{D5.17} created on $\big[\ct^c\big]\op$.
Combining these we have many examples
of exact equivalences
of triangulated categories $\fs\left(\fs(\cs)\op\right)\cong\cs\op$,
which are homeomorphisms with respect to the metrics.
Thus Problem~\ref{P6.1} asks the reader to find the right generalization.
\emxm

Next one can wonder about the functoriality of the construction. Suppose
$F:\cs\la\ct$ is a triangulated functor, and that both $\cs$ and $\ct$
have good metrics. What are reasonable sufficient conditions which 
guarantee the existence of an induced functor $\fs(F)$, either from
$\fs(\cs)$ to $\fs(\ct)$ or in the other direction? So far
there is one known result of this genre, the reader can find
the statement below
in Sun and Zhang~\cite[Theorem~1.1(3)]{Sun-Zhang21}.

\thm{T6.3}
Suppose we are given two triangulated categories $\cs$ and $\ct$, both
with good metrics. Suppose we are also given a pair of functors
$\xymatrix{
F\colon\cs\ar@<0.5ex>[r] &\ct \ar@<0.5ex>[l]\colon G
}$
with $F\dashv G$, meaning that $F$ is 
left adjoint to $G$. Assume further that both
$F$ and $G$ are continuous with respect to the metrics, in the obvious
sense.

Then the functor $\wh F:\MMod\ct\la\MMod\cs$ induced
by composition with $F$, that is the 
functor taking the $\ct$--module $H:\ct\op\la\ab$ to the 
$\cs$--module $(H\circ F):\cs\op\la\ab$,
restricts to a functor which we will denote
$\fs(F):\fs(\ct)\la\fs(\cs)$. That is the functor $\fs(F)$ is defined
to be the unique map making the square below commute
\[\xymatrix{
\fs(\ct)\ar[r]^-{\fs(F)}\ar@{^{(}->}[d] & \fs(\cs)\ar@{^{(}->}[d]\\
\MMod\ct\ar[r]^-{\wh F} & \MMod\cs
}\]
where the vertical inclusions are given by the definition of
$\fs(?)\subset\fl(?)\subset\MMod?$ of Discussion~\ref{D5.17}~(1) and (2).

Furthermore: the functor $\fs(F)$ respects the exact triangles
as defined in  
Discussion~\ref{D5.17}(3).
\ethm

Sun and Zhang go on to study recollements. Suppose we are 
given a recollement
of triangulated categories
\[\xymatrix@C+20pt{
\car\ar[r]|I &
\cs\ar@/^0.5pc/@<0.3ex>[l]^{I_\rho}\ar@/_0.5pc/@<-0.3ex>[l]_{I_\lambda}\ar[r]|J  &
\ct
\ar@/^0.5pc/@<0.3ex>[l]^{J_\rho}\ar@/_0.5pc/@<-0.3ex>[l]_{J_\lambda}
}\]
If all three triangulated categories come with good metrics, and if all
six functors are continuous, then the following 
may be found in~\cite[Theorem~1.2]{Sun-Zhang21}.

\thm{T6.5}
Under the hypotheses above, applying $\fs$ yields a right recollement
\[\xymatrix@C+20pt{
\fs(\car)\ar@/^0.5pc/@<0.3ex>[r]^{\fs(I_\lambda)} &
\fs(\cs)\ar[l]^{\fs(I)}\ar@/^0.5pc/@<0.3ex>[r]^{\fs(J_\lambda)}  &
\fs(\ct)
\ar[l]^{\fs(J)}
}\]
\ethm

In the presence of enough continous adjoints, we deduce that
a semiorthogonal
decomposition of $\cs$ gives rise to a semiorthogonal decomposition
of $\fs(\cs)$. In view of the fact that there are metrics on
$\dperf X$ and $\dcoh(X)$ such that
\[
\fs\big(\dperf X\big)=\dcoh(X),\qquad\fs\left(\dcoh(X)\op\right)={\dperf X}\op
\]
it is natural to wonder how the recent theorem of
Sun and Zhang~\cite[Theorem~1.2]{Sun-Zhang21} compares with the 
older work of Kuznetsov~\cite[Section~2.5]{Kuznetsov08}
and \cite[Section~4]{Kuznetsov11}.

The above shows that, subject to suitable hypotheses, the
construction taking $\cs$ to $\fs(\cs)$ can preserve (some of) the
internal structure on the category $\cs$---for example semiorthogonal
decompositions. This leads naturally to

\plm{P6.7}
What other pieces of the internal structure of $\cs$ are respected
by the construction that passes to $\fs(\cs)$? Under what conditions
are these preserved?
\eplm

Problem~\ref{P6.7} may sound vague, but it can be made precise enough.
For example there is a huge literature dealing with the group of
autoequivalences of the derived categories $\dcoh(X)$. Now as
it happens the metrics for which Remark~\ref{R5.21} gives the 
equivalences
\[
\fs\big(\dperf X\big)\cong\dcoh(X),\qquad
\fs\left(\dcoh(X)\op\right)\cong{\dperf X}\op
\]
can be given (up to equivalence) intrinsic descriptions.
Note that the way we introduced these metrics, in Discussion~\ref{D5.17},
was to use a preferred t-structure on $\ct=\Dqc(X)$ to give on $\ct$ a
metric, unique up to equivalence, and hence
induced metrics on $\ct^c=\dperf X$ and on $\ct^b_c=\dcoh(X)$ which
are also unique up to equivalence. But this description seems to
depend on an embedding into the large category $\ct$.
What I'm asserting now is that there are alternative
descriptions of the same equivalence classes of metrics on
$\ct^c$ and on $\ct^b_c$, which do not use the embedding into $\ct$.
The interested reader can find this in the later sections of
\cite{Neeman18A}.
Anyway: a consequence is that
any autoequivalence, of either $\dperf X$ or of $\dcoh(X)$,
must be continuous with a continuous inverse. Hence the group
of autoequivalences of $\dcoh(X)$ must be isomorphic to the group
of autoequivalences of $\dperf X$. Or more generally: assume $\ct$ is a 
noetherian, approximable triangulated category, where noetherian
has the meaning of \cite[Definition~5.1]{Neeman18A}. Then
the group of exact autoequivalences of $\ct^c$ is canonically
isomorphic to the group of exact autoequivalences of $\ct^b_c$.

Are there similar theorems about t-structures in $\cs$ going to t-structures
in $\fs(\cs)$? Or about stability conditions on $\cs$ mapping to stability 
conditions on $\fs(\cs)$?

We should note that any such theorem will have to come with conditions.
After all: the category $\dcoh(X)$ always has a bounded t-structure,
while Antieau, Gepner and Heller~\cite[Theorem~1.1]{Antieau-Gepner-Heller19}
shows that $\dperf X$ doesn't in general. Thus it is possible
for $\cs$ to have a bounded t-structure but for $\fs(\cs)$ not to.
And in this particular example the equivalence class
of the metric has an intrinsic description, in the sense
mentioned above. 

Perhaps we should remind the reader: the 
article~\cite{Antieau-Gepner-Heller19},
by Antieau, Gepner and Heller, finds a $K$-theoretic obstruction
to the existence of bounded t-structures, more precisely if
an appropriate category $\ce$ has a bounded t-structure then
$K_{-1}(\ce)=0$.
Hence the reference  to~\cite{Antieau-Gepner-Heller19}
immediately  raises
the question of how the construction passing from $\cs$ to $\fs(\cs)$
might relate to $K$-theory, especially to negative $K$-theory. Of 
course: one has to be a little circumspect here. While there is a
$K$-theory for triangulated categories (see \cite{Neeman05A} for a survey),
this $K$-theory has only been proved to behave 
well for ``nice'' triangulated categories,
for example for triangulated categories with bounded t-structures.
Invariants like negative $K$-theory have never been defined for
triangulated categories, and might well give nonsense. In what follows
we will assume that all the $K$-theoretic statements are for triangulated
categories with chosen enhancements, and that $K$-theory means
the Waldhausen $K$-theory of the enhancement. We recall in passing that
the enhancements are unique for many interesting classes
of triangulated categories, see Lunts and Orlov~\cite{Lunts-Orlov10},
Canonaco and Stellari~\cite{Canonaco-Stellari18}, Antieau~\cite{Antieau18}
and Canonaco, Neeman and Stellari~\cite{Canonaco-Neeman-Stellari21}.

With the disclaimers out of the way: what do
the results surveyed in this article have to do with negative $K$-theory?

Let us begin with Schlichting's 
conjecture~\cite[Conjecture~1 of Section~10]{Schlichting06}; this 
conjecture, now known to be false~\cite{Neeman21}, predicted that
the negative $K$-theory of any abelian category should vanish.
But Schlichting also proved that (1) $K_{-1}(\ca)=0$ for
any abelian category $\ca$, and (2) $K_{-n}(\ca)=0$ whenever $\ca$ is
a noetherian abelian category and $n>0$. 
Now note that the $K(\ca)=K(\ca\op)$, hence
the negative $K$-theory of any artinian abelian category must also 
vanish. And playing with extensions of abelian categories, we easily
deduce the vanishing of the negative $K$-theory of a sizeable class
of abelian categories. So while Schlichting's conjecture is false 
in the generality in which it was stated, there is some large class
of abelian categories for which it's true. 
The challenge is to understand this class.

It becomes interesting to see what relation, if any, the results surveyed
here have with this question.

Let us begin with Theorems~\ref{T3.7} and \ref{T3.3}. Theorem~\ref{T3.7}
tells us that, when $X$ is a quasiexcellent, finite-dimensional,
separated noetherian scheme, the category
$\dcoh(X)$ is strongly generated. This category has a
unique enhancement whose $K$-theory agrees with the $K$-theory of
the noetherian abelian category $\text{Coh}(X)$, hence the 
negative $K$-theory vanishes. Theorem~\ref{T3.3}
and Remark~\ref{R3.5} tell us that
the category $\dperf X$ has a strong generator if and only if
$X$ is regular and finite dimensional---in which case it is 
equivalent to $\dcoh(X)$ and its unique enhancement has vanishing
negative $K$-theory. This raises the question:

\plm{P6.9}
If $\ct$ is a triangulated category with a strong
generator, does it follow
that any enhancement of $\ct$ has vanishing negative $K$-theory?
\eplm

Let us refine this question a little. In Definition~\ref{D3.1} we learned
that a \emph{strong generator,} for a triangulated category $\ct$, is an
object $G\in\ct$ such that there exists an integer $\ell>0$ with
$\ct=\genuf G\ell$. Following Rouquier,
we can ask for estimates on the integer $\ell$. This leads
us to:

\dfn{D6.11}
Let $\ct$ be a triangulated category. The \emph{Rouquier dimension} of $\ct$
is the smallest integer $\ell\geq0$ (we allow 
the possibility $\ell=\infty$), for which
there exists an object $G$ with $\ct=\genuf G{\ell+1}$. See
Rouquier~\cite{Rouquier08} for much more about this fascinating invariant.
\edfn 

There is a rich and beautiful literature estimating this invariant
and its various cousins---see Rouquier~\cite{Rouquier08} for the 
origins of the theory, and a host of other 
places for subsequent developments.
For this survey we note only that, for $\dcoh(X)$, the Rouquier
dimension is conjectured to be equal to the Krull dimension of $X$.
But by a conjecture of Weibel~\cite{Weibel80}, now a theorem
of Kerz, Strunk and Tamme~\cite{Kerz-Strunk-Tamme18}, the 
Krull dimension
of $X$ also has a $K$-theoretic description: the groups $K_n$
of the unique enhancement of $\dperf X$ vanish for all $n<-\dim(X)$.
Recalling that $\cs=\dcoh(X)$ is related to $\dperf X$ by the fact that
the construction $\fs$ interchanges them (up to passing to opposite
categories, which has no effect on $K$-theory),
this leads us to ask:

\plm{P6.13}
Let $\cs$ be a regular ($=$ strongly generated) triangulated category as in
Definition~\ref{D3.1}, and let $N<\infty$ be its Rouquier
dimension. Is it true that $K_n$ vanishes on any enhancement
of $\fs(\cs)$, for any metric on $\cs$ and whenever $n<-N$?  
\eplm 

In an entirely different direction: we know that the construction 
$\fs$ interchanges $\dperf X$ and $\dcoh(X)$, and that these categories
coincide if and only if $X$ is regular. This leads us to ask:

\plm{P6.21}
Is there a way to measure the ``distance'' between $\cs$ and $\fs(\cs)$,
in such a way that resolution of singularities can be viewed as a
process reducing this distance? Who knows: maybe there is even a good 
metric on
$\cs=\dperf X$ and/or on $\cs'=\dcoh(X)$, such that the construction
$\fs$ takes either $\cs$ or $\cs'$ to an $\fs(\cs)$ or $\fs(\cs')$
which is $\dperf Y=\dcoh(Y)$ for some resolution of singularities $Y$
of $X$.
\eplm

While on the subject of regularity ($=$ strong generation): 

\plm{P6.22}
Is there some way to understand
which are the approximable triangulated categories $\ct$
for which $\ct^c$ and/or
$\ct^b_c$ are regular? 
\eplm

Theorems~\ref{T3.3} and \ref{T3.7}, deal with the case $\ct=\Dqc(X)$.
Approximability is used in the proofs given
in~\cite{Neeman17} and~\cite{Aoki20}, but only to 
ultimately reduce to
the case of $\ct^c=\dperf X$ with $X$ an affine scheme---this
case turns
out to be classical, it was settled already
in Kelly's 1965 article~\cite{Kelly65}. And the diverse
precursors of
Theorems~\ref{T3.3} and \ref{T3.7}, which we have hardly mentioned
in the current survey, are also relatively narrow in scope. 
But presumably there are other proofs out there, yet to be discovered. 
And new approaches 
might well lead to generalizations that hold for triangulated
categories having nothing to do with algebraic geometry.

Next let's revisit Theorem~\ref{T29.27}, the theorem identifying
each of 
$\left[\ct^c\right]\op$ (respectively $\ct^b_c$) as the finite 
homological functors on the other. In view of the
motivating application, discussed in Remark~\ref{R5.999},
and of the generality of Theorem~\ref{T29.27},
it's natural to wonder: 

\plm{P6.23}
Do GAGA-type theorems have interesting
generalizations to other approximable triangulated categories?
The reader is invited to check \cite[Section~8 and
Appendix~A]{Neeman17A}: except for the couple of paragraphs 
in \cite[Example~A.2]{Neeman17A} everything is formulated in gorgeous
generality and might be applicable in other contexts.

In the context of $\dcoh(X)$, where $X$ is a
scheme proper over a noetherian ring $R$, there
was a wealth of different-looking GAGA-statements before Jack Hall's
lovely 
paper~\cite{Hall18} unified them into one. In other words: 
the category $\dcoh(X)=\ct^b_c$
had many different-looking incarnations, and it wasn't until Hall's
paper that it was understood that there was one underlying reason 
why they all coincided. 

Hence Problem~\ref{P6.23} asks whether this pattern is present for 
other $\ct^b_c$, in other words for $\ct^b_c\subset\ct$ where $\ct$
are some other $R$--linear, approximable triangulated categories.
\eplm

And finally:

\plm{P6.25}
Is there a version of Theorem~\ref{T29.27} that holds for non-noetherian
rings? 
\eplm

There is evidence that something might be true, see
Ben-Zvi, Nadler and
Preygel~\cite[Section~3]{BenZvi-Nadler-Preygel16}.
But
the author has no idea what the right statement ought to be, let
alone how to go about proving it.

%------
% Insert acknowledgments and information
% regarding funding at the end of the last
% section, i.e., right before the bibliography.
%------

\begin{ack}
We would like to thank Asilata Bapat, Jim Borger, 
Anand Deopurkar, Jack Hall, Bernhard Keller, 
Tony Licata and Bregje Pauwels 
for corrections and improvements to earlier drafts. Needless to
say, the flaws that remain are entirely the author's fault.
\end{ack}

\begin{funding}
This work was partially supported by the Australian Research
Council, more specifically by grants number DP150102313 and
DP200102537.
\end{funding}

%\begin{thebibliography}{99}

%------ Example for a paper in journal:
% \bibitem{article1}
% A.~Petrunin, Parallel transportation for Alexandrov space with curvature bounded below.
% \emph{Geom. Funct. Anal.} \textbf{8} (1998), no.~1, 123--148.

%------ Example for a book:
% \bibitem{book1}
% W.~P. Ziemer, \emph{Weakly differentiable functions}.
% Grad. Texts in Math. 120,  Springer, New York, 1989.

%------ Example for a paper in a book:
% \bibitem{incollection1}
% J.~S. Milne, Introduction to Shimura varieties.
% In \emph{Harmonic analysis, the trace formula, and Shimura varieties},
% edited by M.~W. Marcellin and E.~Giorgi, pp. 265--378,
% Clay Math. Proc. 4, Amer. Math. Soc., Providence, RI, 2005.

%------ Example for a preprint on arXiv:
% \bibitem{preprint1}
% D.~V. Nguyen, S.~K. Chilappagari, M.~W. Marcellin, and B.~Vasic,
% LDPC codes from latin squares free of small trapping sets,
% 2010, \href{http://arxiv.org/abs/1008.4177}{arXiv:1008.4177}.

%------ Example for a report:
% \bibitem{report1}
% J.~Schöberl, Commuting quasi-interpolation operators.
% Technical report isc-01-10-math, Texas A\&M University, 2001,
% \url{www.isc.tamu.edu/publications-reports/tr/0110.pdf}.

%------ Example for a thesis:
% \bibitem{thesis1}
% E.~Giorgi, \emph{The geometric universe}.
% Ph.D. thesis, University of Maryland, College Park, 2002.

%\end{thebibliography}

\end{document}